\documentclass{amsart}
\usepackage{amssymb}

\newcommand{\Z}{\mathbb{Z}}

\newcommand{\C}{\mathbb{C}}

\newcommand{\Fq}{\mathbb{F}_q}
\newcommand{\Fqs}{\mathbb{F}_{q^{2}}}
\newcommand{\Fqe}{\mathbb{F}_{q^{e}}}

\newcommand{\CalP}{\mathcal{P}}

\newcommand{\bbP}{\mathbb{P}}

\newcommand{\Ind}{\mathrm{Ind}}

\newcommand{\inv}{\mathrm{inv}}

\newcommand{\bG}{\overline{G}}
\newcommand{\bT}{\overline{T}}
\newcommand{\bK}{\overline{K}}
\newcommand{\bH}{\overline{H}}
\newcommand{\tsigma}{{\tilde{\sigma}}}

\newcommand{\unnu}{{\underline{\nu}}}

\newcommand{\unrho}{{\underline{\rho}}}

\newcommand{\isomto}{\overset{\sim}{\rightarrow}}
\newcommand{\leftdiv}{\!\setminus\!}
\newcommand{\sigL}{\langle\sigma\rangle\leftdiv L}
\newcommand{\subsigL}{\langle\sigma\rangle\setminus L}
\newcommand{\ffi}{\mathrm{ff-inv}}
\newcommand{\lng}{\langle}
\newcommand{\rng}{\rangle}

\numberwithin{equation}{section}

\newtheorem{theorem}{Theorem}[section]
\newtheorem{corollary}[theorem]{Corollary}
\newtheorem{lemma}[theorem]{Lemma}
\newtheorem{proposition}[theorem]{Proposition}

\theoremstyle{definition}
\newtheorem*{definition}{Definition}
\newtheorem*{example}{Example}

\newtheorem*{ack}{Acknowledgements}

\newcommand{\bdf}{\begin{definition}}
\newcommand{\edf}{\end{definition}\noindent}
\newcommand{\bex}{\begin{example}}
\newcommand{\eex}{\end{example}\noindent}
\newcommand{\blm}{\begin{lemma}}
\newcommand{\elm}{\end{lemma}}
\newcommand{\bpr}{\begin{proposition}}
\newcommand{\epr}{\end{proposition}}
\newcommand{\bth}{\begin{theorem}}
\renewcommand{\eth}{\end{theorem}}
\newcommand{\bpf}{\begin{proof}}
\newcommand{\epf}{\end{proof}\noindent}
\newcommand{\bcr}{\begin{corollary}}
\newcommand{\ecr}{\end{corollary}\noindent}
\newcommand{\beq}{\begin{equation}}
\newcommand{\eeq}{\end{equation}}
\newcommand{\bes}{\begin{equation*}}
\newcommand{\ees}{\end{equation*}}
\newcommand{\ben}{\begin{enumerate}}
\newcommand{\een}{\end{enumerate}}

\begin{document}
\title[Induced characters of the projective general linear group]
{Induced characters of the projective general linear group
over a finite field}

\author{Anthony Henderson}
\address{School of Mathematics and Statistics,
University of Sydney, NSW 2006, AUSTRALIA}
\email{anthonyh@maths.usyd.edu.au}
\thanks{This work was supported by Australian Research Council grant DP0344185}
\begin{abstract}
Using a general result of Lusztig, we find the decomposition into
irreducibles of certain induced characters of the projective general
linear group over a finite field of odd characteristic.
\end{abstract}
\maketitle
\section{Introduction}
Let $q$ be a power of an \textbf{odd} prime, and $\Fq$ the finite field
with $q$ elements. Fix an \textbf{even} integer $n\geq 2$. Consider
the ordinary (complex) character theory of the projective general linear group
$PGL_n(\Fq)=\{\bar{g}\,|\,g\in GL_n(\Fq)\}$, where $\bar{g}$ denotes the coset
of $g$ relative to $Z(GL_n(\Fq))\cong\Fq^\times$.
The aim of this paper is
to calculate the decomposition into irreducible characters of the following
induced characters:
\[ \Ind_{PGSp_n(\Fq)}^{PGL_n(\Fq)}(1),\ 
\Ind_{PGO_n^+(\Fq)}^{PGL_n(\Fq)}(1),\text{ and }
\Ind_{PGO_n^-(\Fq)}^{PGL_n(\Fq)}(1). \]
Specifically, for each of these induced characters $I$ and for each
irreducible character $\chi$ we will give
a (manifestly nonnegative) combinatorial formula for the multiplicity
$\langle\chi,I\rangle$, in the same spirit as the formulas given in
\cite{mysymm} for the corresponding induced characters of $GL_n(\Fq)$.

Recall the definitions of these subgroups of $PGL_n(\Fq)$ (which are
in fact only defined up to conjugacy). The group $PGL_n(\Fq)$ acts on the
set
\[ PGL_n(\Fq)^{\mathrm{symm}}=\{\bar{h}\in PGL_n(\Fq)\,|\,\overline{h^{t}}
=\bar{h}\}
=\{\bar{h}\,|\, h\in GL_n(\Fq), h^{t}
=\pm h\} \]
by the rule $\bar{g}.\bar{h}=\overline{ghg^t}$; we can think of 
$PGL_n(\Fq)^{\mathrm{symm}}$
as the set of equivalence classes of nondegenerate
symmetric and skew-symmetric bilinear forms on $\Fq^n$, where two
forms are equivalent if they are scalar multiples of each other.
The subgroup $PGSp_n(\Fq)$ (also called $PCSp_n(\Fq)$) is the
stabilizer of an equivalence class of nondegenerate skew-symmetric forms.
The subgroup $PGO_n^+(\Fq)$ is the stabilizer of an equivalence class
of nondegenerate symmetric forms whose Witt index is $n/2$, and
$PGO_n^-(\Fq)$ is the stabilizer of an equivalence class
of nondegenerate symmetric forms whose Witt index is $n/2-1$.
Since these three types of forms precisely give the three orbits
of $PGL_n(\Fq)$ on $PGL_n(\Fq)^{\mathrm{symm}}$, the character of the
corresponding permutation representation is just the sum of the three induced
characters we are calculating.

The reason for excluding the case when $n$ is odd is that it is
already known: $PGL_n(\Fq)$ acts transitively on $PGL_n(\Fq)^{\mathrm{symm}}$,
and the stabilizer $PGO_n(\Fq)$ is the image of 
$O_n(\Fq)=SO_n(\Fq)\times\{\pm 1\}$; in other words,
$PGO_n(\Fq)=PSO_n(\Fq)$. Hence if $n$ is odd,
\beq 
\langle\chi,\Ind_{PGO_n(\Fq)}^{PGL_n(\Fq)}(1)\rangle
=\langle\chi,\Ind_{SO_n(\Fq)}^{GL_n(\Fq)}(1)\rangle
=\langle\chi,\Ind_{O_n(\Fq)}^{GL_n(\Fq)}(1)\rangle,
\eeq
where $\chi$ is regarded as an irreducible character of $GL_n(\Fq)$
via the canonical projection $GL_n(\Fq)\to PGL_n(\Fq)$. Thus
\cite[Theorem 4.1.1]{mysymm} gives a formula
for this multiplicity.

By contrast, when $n$ is even, the image in $PGL_n(\Fq)$ of the stabilizer of
$h$ in $GL_n(\Fq)$, namely
\[ \{\bar{g}\in PGL_n(\Fq)\,|\,ghg^{t}=h\}
=\{\bar{g}\in PGL_n(\Fq)\,|\,ghg^{t}=\alpha h,\,
\exists\alpha\in(\Fq^\times)^2\}, \]
is a subgroup of index $2$ in the stabilizer of $\bar{h}$ in $PGL_n(\Fq)$,
namely
\[ \{\bar{g}\in PGL_n(\Fq)\,|\,ghg^t=\alpha h,\,\exists\alpha\in\Fq^\times\}.\]
(The fact that $\alpha$ can take non-square values can be seen directly for
$n=2$ and then deduced for general even $n$.)
If $\omega$ denotes, in each case, the homomorphism from the latter 
stabilizer to $\{\pm 1\}$ whose kernel is this subgroup, we have
\beq
\begin{split}
\langle\chi,\Ind_{PGSp_n(\Fq)}^{PGL_n(\Fq)}(1)\rangle
+\langle\chi,\Ind_{PGSp_n(\Fq)}^{PGL_n(\Fq)}(\omega)\rangle
&=\langle\chi,\Ind_{Sp_n(\Fq)}^{GL_n(\Fq)}(1)\rangle,\\
\langle\chi,\Ind_{PGO_n^+(\Fq)}^{PGL_n(\Fq)}(1)\rangle
+\langle\chi,\Ind_{PGO_n^+(\Fq)}^{PGL_n(\Fq)}(\omega)\rangle
&=\langle\chi,\Ind_{O_n^+(\Fq)}^{GL_n(\Fq)}(1)\rangle,\\
\langle\chi,\Ind_{PGO_n^-(\Fq)}^{PGL_n(\Fq)}(1)\rangle
+\langle\chi,\Ind_{PGO_n^-(\Fq)}^{PGL_n(\Fq)}(\omega)\rangle
&=\langle\chi,\Ind_{O_n^-(\Fq)}^{GL_n(\Fq)}(1)\rangle.
\end{split}
\eeq
So the formulas we are seeking cannot be deduced from the results 
in \cite{mysymm},
but once found they can be combined with those results to give 
formulas for the induced characters
$\Ind(\omega)$ in the three cases.
(If $q\equiv 1$ mod $4$ and $n\equiv 2$ mod $4$, this
gives nothing substantially different, because $\omega$ is then the restriction
of the unique nontrivial homomorphism $PGL_n(\Fq)\to\{\pm 1\}$.)

Since $\Ind_{Sp_n(\Fq)}^{GL_n(\Fq)}(1)$ is multiplicity-free,
$\Ind_{PGSp_n(\Fq)}^{PGL_n(\Fq)}(1)$ must be also, so our formula merely
specifies which constituents of the former
are consituents of the latter; in the other cases,
we have the inequalities
\beq 
0\leq \langle\chi,\Ind_{PGO_n^\epsilon(\Fq)}^{PGL_n(\Fq)}(1)\rangle
\leq \langle\chi,\Ind_{O_n^\epsilon(\Fq)}^{GL_n(\Fq)}(1)\rangle,
\eeq
which will be obviously satisfied by our formulas.

Our approach will be the same as that of \cite{mysymm}: namely, we use Lusztig's
formula \cite[Theorem 3.3]{symmfinite} for the inner product
$\lng R_{\bT}^{\lambda},I\rng$, where $R_{\bT}^{\lambda}$ is a Deligne-Lusztig
character of $PGL_n(\Fq)$ and $I$ is one of the above induced characters.
Since the irreducible characters are linear combinations of the Deligne-Lusztig
characters with known coefficients, all that remains is to interpret Lusztig's
formula combinatorially and to manipulate it using various facts about characters
of symmetric groups.

This paper does not attempt to be self-contained, and reference is freely
made to the arguments and definitions of \cite{mysymm}. However, the more
limited scope of the present inquiry (no unitary groups, no inner involutions)
enables slight simplifications in the exposition and in the notation, which
I hope will make the reader's job easier.

Since the complete formulas (Theorems \ref{firstpgspthm} and \ref{firstpgothm} 
below) cannot be stated until we have introduced further notation,
here are the answers in the
special case of \textbf{unipotent} characters. Recall that the unipotent
characters of $GL_n(\Fq)$, all of which factor through $PGL_n(\Fq)$,
are parametrized by partitions of $n$, $\rho\mapsto\chi^\rho$; in our
conventions, $\chi^{(n)}$ is the trivial character and $\chi^{(1^n)}$ is
the Steinberg character. We have
\beq
\langle\chi^\rho,\Ind_{PGSp_n(\Fq)}^{PGL_n(\Fq)}(1)\rangle
=\langle\chi^\rho,\Ind_{Sp_n(\Fq)}^{GL_n(\Fq)}(1)\rangle
=\left\{\begin{array}{cl}
1,&\text{ if $\rho$ is even}\\
0,&\text{ otherwise.}
\end{array}\right. 
\eeq
In other words, the unipotent constituents of
$\Ind_{Sp_n(\Fq)}^{GL_n(\Fq)}(1)$ are all constituents of
$\Ind_{PGSp_n(\Fq)}^{PGL_n(\Fq)}(1)$ also.
As for the other cases, recall from \cite{mysymm} that
\beq
\langle\chi^\rho,\Ind_{O_n^\epsilon(\Fq)}^{GL_n(\Fq)}(1)\rangle
=\frac{1}{2}\prod_i (m_i(\rho)+1) +\left\{\begin{array}{cl}
\frac{1}{2}\epsilon,&\text{ if $\rho'$ is even}\\
0,&\text{ otherwise.}
\end{array}\right.
\eeq
Here $m_i(\rho)$ denotes the multiplicity of $i$ as a part of $\rho$,
and $\rho'$ is the transpose partition (which is even if and only if
$m_i(\rho)$ is even for all $i$). The new formula is:
\beq
\begin{split}
\langle\chi^\rho,\Ind_{PGO_n^\epsilon(\Fq)}^{PGL_n(\Fq)}(1)\rangle
&=\frac{1}{4}\prod_i (m_i(\rho)+1) +\left\{\begin{array}{cl}
\frac{1}{2}\epsilon,&\text{ if $\rho'$ is even}\\
0,&\text{ otherwise}
\end{array}\right.\\
&\quad+\left\{\begin{array}{cl}
\frac{(-1)^{\frac{1}{2}\ell(\rho)_1}}{4}\prod_i(m_{2i}(\rho)+1),
&\text{ if $2\mid m_{2i+1}(\rho)$, $\forall i$}\\
0,&\text{ otherwise.}
\end{array}\right.
\end{split}
\eeq
Here $\ell(\rho)_1=\sum_i m_{2i+1}(\rho)$ is the number of odd parts of $\rho$.
\bex
The multiplicities of the unipotent characters of $PGL_4(\Fq)$ in the
above induced characters are given in the following table:
\[ \begin{array}{c|ccccc|}
 & (4) & (31) & (2^2) & (21^2) & (1^4)\\
\hline
\rule[-0.3cm]{0cm}{0.6cm}
\Ind_{PGSp_4(\Fq)}^{PGL_4(\Fq)}(1) & 1 & 0 & 1 & 0 & 0\\ 
\rule[-0.3cm]{0cm}{0.6cm}
\Ind_{PGO_4^+(\Fq)}^{PGL_4(\Fq)}(1) & 1 & 1 & 2 & 1 & 2\\
\rule[-0.3cm]{0cm}{0.6cm}
\Ind_{PGO_4^-(\Fq)}^{PGL_4(\Fq)}(1) & 1 & 1 & 1 & 1 & 1\\
\hline
\end{array} \]
See \eqref{exampleeqn} below for the complete decomposition in the case $n=2$.
\eex

\begin{ack}
The stimulus for this work was a question of Jorge Soto-Andrade, and
I would like to thank him and Dipendra Prasad for their interest.
\end{ack}
\section{Notation and Statement of Results}
For general combinatorial notation, we follow \cite{macdonald}.
For instance, we write $\mu\vdash m$ to mean that $\mu$ is a partition
of $m$; the size $m$ can then be written $|\mu|$. The nonzero
parts of $\mu$ are $\mu_1\geq\mu_2\geq\cdots\geq
\mu_{\ell(\mu)}$. The transpose partition is
$\mu'$. For $i\in\Z^+$,
the multiplicity of $i$ as a part of $\mu$ is written
$m_i(\mu)$. We say that $\mu$ is even if all its parts
are even, or equivalently if $2\mid m_{i}(\mu'),\ \forall i$.

We write $w_\mu$ for a permutation in $S_{|\mu|}$ whose cycle-type is $\mu$
(determined up to conjugacy), $\epsilon_\mu\in\{\pm 1\}$ for its sign,
and $z_\mu$ for the order of its centralizer in $S_{|\mu|}$.

As well as $\ell(\mu)$ for the length of $\mu$, we will use the notations
$\ell(\mu)_0$ and $\ell(\mu)_1$ for the number of even and odd
(nonzero) parts of $\mu$ respectively. 
Note that $\epsilon_\mu=(-1)^{\ell(\mu)_0}$.
On occasion we will need to further analyse $\ell(\mu)_0$
into $\ell(\mu)_{0\,\mathrm{mod}\,4}$ and $\ell(\mu)_{2\,\mathrm{mod}\,4}$.

For any finite group $\Gamma$, $\widehat{\Gamma}$ denotes the set of
irreducible (complex) characters of $\Gamma$.
We have $\widehat{S_m}=\{\chi^{\rho}\,|\,\rho\vdash m\}$, where
$\chi^\rho$ is as in \cite[I.7]{macdonald}.
Write $\chi_\mu^{\rho}$ for the value of $\chi^{\rho}$ at an element
of cycle-type $\mu$,
so that $\chi_\mu^{(m)}=1$ and 
$\chi_\mu^{\rho'}=\epsilon_\mu\chi_\mu^{\rho}$.

Throughout the paper, we fix an algebraic closure $k$ of $\Fq$, 
a $k$-vector space $V$ of dimension $n$ (our fixed \textbf{even} positive 
integer) and a Frobenius map $F:V\to V$ relative to $\Fq$.
Let $G=GL(V)$, $\bG=PGL(V)$, and write $F$ also for the induced
Frobenius maps on these groups. Thus $G^F\cong GL_n(\Fq)$ and
$\bG^F\cong PGL_n(\Fq)$. For any subgroup $H$ of $G$, 
we write $\overline{H}$ for the image
of $H$ under the canonical projection
$G\to \bG$; in other words, $\overline{H}=HZ/Z$ where $Z=Z(G)\cong k^\times$
is the kernel. Because the Lang map of $Z$
is surjective, we have a short exact sequence
$1\to Z^F\to (HZ)^F\to \bH^F\to 1$.
(But note that $(HZ)^F$ is in general larger than $H^F Z^F$, so
$H^F\to\bH^F$ is not necessarily surjective.)

As in \cite[Chapter IV]{macdonald},
we need to define the `dual group' of $k^\times$: consider the
system of maps $\widehat{\Fqe^{\times}}\to
\widehat{\mathbb{F}_{q^{e'}}^{\times}}$ for $e\,|\,e'$
(the transpose of the norm map), and its direct limit
$L=\underset{\longrightarrow}{\mathrm{lim}}\ \widehat{\Fqe^{\times}}$.
Let $\sigma$ denote the $q$-th power map
on both $k^{\times}$ and $L$, so that 
$(k^{\times})^{\sigma^{e}}\cong\Fqe^{\times}$,
$L^{\sigma^{e}}\cong\widehat{\Fqe^{\times}}$ for all $e\geq 1$.
Write $\langle\cdot,\cdot\rangle_e:
\Fqe^{\times}\times L^{\sigma^{e}} \to \C^{\times}$
for the canonical pairing, which is `bimultiplicative' in the sense that
\[ \lng\alpha_1\alpha_2,\xi\rng_e=\lng\alpha_1,\xi\rng_e\lng\alpha_2,\xi\rng_e,\
\lng\alpha,\xi_1\xi_2\rng_e=\lng\alpha,\xi_1\rng_e\lng\alpha,\xi_2\rng_e, \]
for all $\alpha,\alpha_1,\alpha_2\in\Fqe^{\times}$, $\xi,\xi_1,\xi_2\in L^{\sigma^e}$. 
We will use without comment the following rules relating the pairings
$\langle\cdot,\cdot\rangle_d$ and $\langle\cdot,\cdot\rangle_e$ when $d\mid e$:
\begin{itemize}
\item if $\alpha\in\mathbb{F}_{q^{d}}^{\times}$ and $\xi\in L^{\sigma^e}$, then
\[ \lng\alpha,\xi\rng_e=\lng\alpha,\xi^{1+q^d+q^{2d}+\cdots+q^{e-d}}\rng_d; \]
\item if $\alpha\in\mathbb{F}_{q^{e}}^{\times}$ and $\xi\in L^{\sigma^d}$, then
\[ \lng\alpha,\xi\rng_e=\lng\alpha^{1+q^d+q^{2d}+\cdots+q^{e-d}},\xi\rng_d; \]
\item if $\alpha\in\mathbb{F}_{q^{d}}^{\times}$ and $\xi\in L^{\sigma^d}$, then
\[ \lng\alpha,\xi\rng_e=\lng\alpha,\xi\rng_d^{e/d}. \]
\end{itemize}
When $e=1$ we omit the subscript in $\lng\cdot,\cdot\rng_e$
and write simply $\langle\cdot,\cdot\rangle$.

We fix some set of representatives for the orbits of the group
$\langle\sigma\rangle$ generated by $\sigma$ on $L$, and call it
$\langle\sigma\rangle\!\setminus\! L$. For
$\xi\in\langle\sigma\rangle\!\setminus\! L$, let
$m_\xi=|\langle\sigma\rangle.\xi|$, in other words the smallest $e\geq 1$
such that $\xi^{q^e}=\xi$. Let 
$N(\xi)=\xi^{1+q+q^2+\cdots+q^{m_\xi-1}}\in L^\sigma$.
Let $d_\xi=\lng -1,N(\xi)\rng=
\langle -1,\xi \rangle_{m_\xi}$, which is
$1$ if $L^{\sigma^{m_\xi}}$ contains square roots of $\xi$, and $-1$ if it
does not. (Calling this sign $d_\xi$ is the piece of notation in \cite{mysymm}
I most regret, but it would be too confusing to change now.) Let
$\eta$ denote the unique element of $L^\sigma$ with order $2$; note that
$d_\eta=1$ if and only if $q\equiv 1$ (mod $4$).

Let $\widehat{\CalP}_n$ be the set of collections
of partitions $\underline{\nu}=(\nu_\xi)_{\xi\in L}$, almost all zero,
such that $\sum_{\xi\in L} |\nu_\xi|= n$.
Let $\widehat{\CalP}_n^{\sigma}$ be the subset of $\widehat{\CalP}_n$
of all $\underline{\nu}$ such that $\nu_{\sigma(\xi)}=\nu_\xi$ for all $\xi$.
Note that for $\underline{\nu}\in\widehat{\CalP}_n^{\sigma}$,
\beq 
\sum_{\xi\in\langle\sigma\rangle\setminus L} m_\xi |\nu_\xi| = n.
\eeq
For $\underline{\nu},\underline{\rho}\in\widehat{\CalP}_n^{\sigma}$,
we write $|\underline{\nu}|=|\underline{\rho}|$ to mean that
$|\nu_\xi|=|\rho_\xi|$ for all $\xi$. For
$\unnu\in\widehat{\CalP}_n^{\sigma}$ we define
\[ \Pi(\unnu):=\prod_{\xi\in L}\xi^{|\nu_\xi|}
=\prod_{\xi\in\subsigL}N(\xi)^{|\nu_\xi|}\in L^\sigma, \]
and let $\widehat{\bbP}_n^{\sigma}$ denote the subset of 
$\widehat{\CalP}_n^\sigma$
consisting of all $\unnu$ such that $\Pi(\unnu)=1$. Finally, if
$\unnu\in\widehat{\bbP}_n^{\sigma}$ is such that $m_\xi|\nu_\xi|$ is always even,
we define
\[ \Phi(\unnu):=\prod_{\xi\in\subsigL} 
\lng\sqrt{\beta},\xi\rng_{m_\xi|\nu_\xi|} \]
for some $\sqrt{\beta}\in\Fqs\setminus\Fq$ whose square $\beta$ lies in $\Fq$.
We have $\Phi(\unnu)\in\{\pm 1\}$ since
its square is $\lng\beta,\Pi(\unnu)\rng=1$. Moreover, $\Phi(\unnu)$
is independent of the element $\sqrt{\beta}$ used to define it, 
because multiplying $\sqrt{\beta}$
by $\gamma\in\Fq^\times$ multiplies $\Phi(\unnu)$ by $\lng\gamma,\Pi(\unnu)\rng=1$.

Now there is a well-known bijection
between $\widehat{\CalP}_n^{\sigma}$ and the set of $G^{F}$-orbits of pairs
$(T,\lambda)$ where $T$ is an $F$-stable maximal torus of $G$
and $\lambda\in\widehat{T^F}$. 
If $(T,\lambda)$ is in the orbit
corresponding to $\underline{\nu}$ under this bijection, we have
the following properties:
\begin{itemize}
\item the fixed lines of $T$ can be labelled
\beq \label{toruseqn} 
\{L_{(\xi,j,i)}\,|\, \xi\in\langle\sigma\rangle\!\setminus\! L,\
1\leq j \leq \ell(\nu_\xi),\ i\in\Z/m_\xi(\nu_\xi)_j\Z\}
\eeq
so that if $t\in T$ has eigenvalue $\alpha_{(\xi,j,i)}$ on each $L_{(\xi,j,i)}$,
$F(t)$ has eigenvalue $\alpha_{(\xi,j,i-1)}^q$ on each $L_{(\xi,j,i)}$;
\item  consequently,
\beq \label{ftoruseqn} 
T^{F} \cong \prod_{\xi\in\langle\sigma\rangle\setminus L}
\prod_{j=1}^{\ell(\nu_\xi)}
\mathbb{F}_{q^{m_\xi(\nu_\xi)_j}}^\times
\eeq
via the map sending $t$ to the collection $(\alpha_{(\xi,j,1)})$;
\item  for $t\in T^F$ as above,
\beq \label{lambdaeqn}
\lambda(t)=
\prod_{\xi\in\langle\sigma\rangle\setminus L}
\prod_{j=1}^{\ell(\nu_\xi)} 
\langle\alpha_{(\xi,j,1)},\xi\rangle_{m_\xi(\nu_\xi)_j}.
\eeq
\end{itemize}
From \eqref{lambdaeqn} it is clear that $\lambda$ is trivial on
$Z^F$ if and only if $\Pi(\unnu)=1$. Hence we obtain a bijection
between $\widehat{\bbP}_n^{\sigma}$ and the set of $\bG^F$-orbits of
pairs $(\bT,\lambda)$ where $\bT$ is an $F$-stable maximal torus of $\bG$
and $\lambda\in\widehat{\bT^F}$. 

For $\underline{\nu}\in\widehat{\CalP}_n^{\sigma}$, let
$B_{\underline{\nu}}$ be the corresponding \textbf{basic character} of 
$G^F$, defined by Green in \cite{green}. 
As proved by Lusztig in \cite{greenpolys},
this coincides with the character of the Deligne-Lusztig virtual
representation $R_{T}^{\lambda}$ for
$(T,\lambda)$ in the corresponding $G^{F}$-orbit. From either point of 
view it is clear that $B_{\unnu}(zg)=\langle z,\Pi(\unnu)\rangle 
B_{\unnu}(g)$
for all $g\in G^F$ and $z\in Z^F$. For 
$\unnu\in\widehat{\bbP}_n^{\sigma}$, $B_\unnu$ may be viewed as a character
of $\bG^F$; it is again the character of the corresponding
Deligne-Lusztig representation $R_{\bT}^\lambda$.

Green's main result on the character theory of $GL_n(\Fq)$ states that
for any $\underline{\rho}\in\widehat{\CalP}_n^{\sigma}$,
\beq \label{mainglneqn}
\chi^{\underline{\rho}} := 
(-1)^{n+\sum_{\xi\in\langle\sigma\rangle\setminus L} |\rho_\xi|}
\sum_{\substack{\underline{\nu}\in\widehat{\CalP}_n^{\sigma}\\
|\underline{\nu}|=|\underline{\rho}|}}
(\prod_{\xi\in\langle\sigma\rangle\setminus L}
(z_{\nu_\xi})^{-1}\chi_{\nu_\xi}^{\rho_\xi}) B_{\underline{\nu}}
\eeq
is an irreducible character of $G^{F}$, and all irreducible
characters arise in this way for unique
$\underline{\rho}\in\widehat{\CalP}_n^{\sigma}$. (Note that our
parametrization differs from that in \cite[Chapter IV]{macdonald}
by transposing all partitions.)
Inverting the transition matrix, we have that
for any $\underline{\nu}\in\widehat{\CalP}_n^{\sigma}$,
\begin{equation} \label{glneqn}
B_{\underline{\nu}}=
(-1)^{n+\sum_{\xi\in\langle\sigma\rangle\setminus L} |\nu_\xi|}
\sum_{\substack{\underline{\rho}\in\widehat{\CalP}_n^{\sigma}\\
|\underline{\rho}|=|\underline{\nu}|}}
(\prod_{\xi\in\langle\sigma\rangle\setminus L}
\chi_{\nu_\xi}^{\rho_\xi}) \chi^{\underline{\rho}}.
\end{equation}
Moreover, $\chi^\unrho(zg)=\langle z,\Pi(\unrho)\rangle \chi^\unrho(g)$
for all $g\in G^F$ and $z\in Z^F$, so the characters which descend
to irreducible characters of $\bG^F$
are exactly $\{\chi^\unrho\,|\,\unrho\in\widehat{\bbP}_n^{\sigma}\}$.
This includes 
the unipotent irreducible characters referred to in the
introduction, which are those $\chi^{\underline{\rho}}$ for which
$\rho_\xi=0$ unless $\xi=1$. (In the introduction we parametrized these
by $\rho=\rho_1$.)

With this notation we can state the results of this paper, which should
be compared with \cite[Theorem 2.1.1]{mysymm} and \cite[Theorem 4.2.1]{mysymm}
respectively.
\bth\label{firstpgspthm}
For any $\unrho\in\widehat{\bbP}_n^{\sigma}$,
\[ \langle \chi^{\underline{\rho}},
\Ind_{PGSp_n(\Fq)}^{PGL_n(\Fq)}(1)\rangle =
\left\{\begin{array}{cl}
1,&\text{ if all $\rho_\xi$ are even and 
$\prod_{\xi\in\langle\sigma\rangle\setminus L} 
N(\xi)^{\frac{|\rho_\xi|}{2}}=1$,}\\
0, &\text{ otherwise.}
\end{array}\right. \]
\eth
\noindent
(Note that if all $\rho_\xi$ are even, 
$\prod_{\xi\in\langle\sigma\rangle\setminus L} 
N(\xi)^{\frac{|\rho_\xi|}{2}}$ is either $1$ or $\eta$,
since its square is $\Pi(\unrho)=1$.)
\bth\label{firstpgothm}
For any $\unrho\in\widehat{\bbP}_n^{\sigma}$ and $\epsilon\in\{\pm 1\}$,
\bes
\begin{split}
\langle \chi^{\underline{\rho}},
\Ind_{PGO_n^\epsilon(\Fq)}^{PGL_n(\Fq)}(1)\rangle &=
\left\{\begin{array}{cl}
{\displaystyle\frac{1}{4}
\prod_{\substack{\xi\in\langle\sigma\rangle\setminus L\\
d_\xi=1}} (\prod_i (m_i(\rho_\xi)+1))},
&\text{ if $d_\xi=-1\Rightarrow\rho_\xi'$ is even}\\
0, &\text{ otherwise}
\end{array}\right.\\
&\thickspace+\left\{\begin{array}{cl}
\frac{1}{2}\epsilon,
&\text{ if all $\rho_\xi'$ are even and 
$\prod_{\xi\in\langle\sigma\rangle\setminus L} 
N(\xi)^{\frac{|\rho_\xi|}{2}}=1$}\\
0, &\text{ otherwise}
\end{array}\right.\\
&\thickspace+\left\{
\begin{array}{rl}
{\displaystyle\frac{\pm 1}{4}
\prod_{\substack{\xi\in\langle\sigma\rangle\setminus L\\
d_\xi=1\\2\nmid m_\xi}}}
&\negthickspace\negthickspace\negthickspace
{\displaystyle(\prod_i(m_{2i}(\rho_\xi)+1))
\prod_{\substack{\xi\in\langle\sigma\rangle\setminus L\\
d_\xi=1\\2\mid m_\xi}}
(\prod_i(m_{i}(\rho_\xi)+1))}\\
&\begin{array}{l}
\text{if $d_\xi=1$,
$2\nmid m_\xi \Rightarrow 2\mid m_{2i+1}(\rho_\xi), \forall i$,}\\
\text{and $d_\xi=-1 \Rightarrow \rho_\xi'$ even}
\end{array}\\ 
0, &\text{ otherwise,} 
\end{array} \right.
\end{split}
\ees
where the sign $\pm 1$ in the third term is
\[ (-1)^{\frac{n}{2}}\,\Phi(\unrho)
\prod_{\substack{\xi\in\langle\sigma\rangle\setminus L\\
d_\xi=1\\2\nmid m_\xi}}(-1)^{\ell(\rho_\xi)_{2\text{ mod }4}}. \]
\eth
\noindent
(Note that the third term can only be nonzero if all $m_\xi|\rho_\xi|$ are even,
so $\Phi(\unrho)$ is defined.)
\bex
The elements of $\widehat{\bbP}_2^{\sigma}$ can be written in `exponential notation'
as follows:
\[ 1^{(2)},\ 1^{(1^2)},\ \eta^{(2)},\ \eta^{(1^2)},\text{ and }
\xi^{(1)}(\xi^{-1})^{(1)}\text{ for
$\xi\in (L^\sigma\cup L^{\tsigma})\setminus\{1,\eta\}$.} \]
Here $L^\tsigma=\{\xi\in L\,|\,\xi^q=\xi^{-1}\}$.
Theorem \ref{firstpgspthm} states that the only irreducible constituent of
$\Ind_{PGSp_2(\Fq)}^{PGL_2(\Fq)}(1)$ is the trivial character $\chi^{1^{(2)}}$,
which is correct because $PGSp_2(\Fq)=PGL_2(\Fq)$. In the following table, each row
contains, for a particular element of $\widehat{\bbP}_2^{\sigma}$, the value of the 
three terms of the right-hand side of Theorem \ref{firstpgothm}, and their sum:
\[ \begin{array}{|c|ccc|c|}
\hline
\rule[-0.6cm]{0cm}{1.2cm}
1^{(2)} & \frac{1}{2} & 0 & \frac{1}{2} & 1\\
\hline
\rule[-0.6cm]{0cm}{1.2cm}
1^{(1^2)} & \frac{3}{4} & \frac{1}{2}\epsilon & -\frac{1}{4} & 
\begin{array}{cl}1,&\text{if }\epsilon=1\\0,&\text{if }\epsilon=-1\end{array}\\
\hline
\rule[-0.6cm]{0cm}{1.2cm}
\eta^{(2)} &
\left\{\begin{array}{cl}\frac{1}{2},&\text{if }d_\eta=1\\0,
&\text{if }d_\eta=-1\end{array}\right\}
& 0 &
\left\{\begin{array}{c}-\frac{1}{2}\\0\end{array}\right\}
& 0\\
\hline
\rule[-0.6cm]{0cm}{1.2cm}
\eta^{(1^2)} &
\left\{\begin{array}{cl}\frac{3}{4},&\text{if }d_\eta=1\\
\frac{1}{4},&\text{if }d_\eta=-1\end{array}\right\}
& 0 &
\left\{\begin{array}{c}\frac{1}{4}\\-\frac{1}{4}\end{array}\right\}
&\begin{array}{cl}1,&\text{if }d_\eta=1\\0,
&\text{if }d_\eta=-1\end{array}\\
\hline
\rule[-0.6cm]{0cm}{1.2cm}
\begin{array}{c}\xi^{(1)}(\xi^{-1})^{(1)},\\
\xi\in L^\sigma\setminus\{1,\eta\}\end{array}
& \left\{\begin{array}{cl}1,&\text{if }d_\xi=1\\0,
&\text{if }d_\xi=-1\end{array}\right\}
& 0 & 0 & \begin{array}{cl}1,&\text{if }d_\xi=1\\0,
&\text{if }d_\xi=-1\end{array}\\
\hline
\rule[-0.6cm]{0cm}{1.2cm}
\begin{array}{c}\xi^{(1)}(\xi^{-1})^{(1)},\\
\xi\in L^{\tsigma}\setminus\{1,\eta\}\end{array}
& \frac{1}{2} & 0 & -\frac{1}{2}\tilde{d}_\xi &
\begin{array}{cl}0,&\text{if }\tilde{d}_\xi=1\\1,
&\text{if }\tilde{d}_\xi=-1\end{array}\\
\hline
\end{array} \]
In the third and fourth rows we have used the fact that $\lng\sqrt{\beta},\eta\rng_2
=\lng-\beta,\eta\rng=-d_\eta$, for $\sqrt{\beta}$ as above.
In the last row, we have used the facts that $d_\xi=1$ and
$\lng\sqrt{\beta},\xi\rng_2=\tilde{d}_\xi$, where $\tilde{d}_\xi$
is the sign in the equation $\sqrt{\xi}^q=\pm\sqrt{\xi}^{-1}$. 
So Theorem \ref{firstpgothm} implies:
\beq \label{exampleeqn}
\begin{split}
\Ind_{PGO_2^+(\Fq)}^{PGL_2(\Fq)}(1)&=\chi^{1^{(2)}}+\chi^{1^{(1^2)}}
+\left\{\begin{array}{cl}\chi^{\eta^{(1^2)}},&\text{if }d_\eta=1\\0,
&\text{if }d_\eta=-1\end{array}\right\}\\
&\qquad\qquad+\sum_{\substack{\xi\in L^\sigma\setminus\{1,\eta\}\\d_\xi=1}}
\chi^{\xi^{(1)}(\xi^{-1})^{(1)}}
+\sum_{\substack{\xi\in L^\tsigma\setminus\{1,\eta\}\\\tilde{d}_\xi=-1}}
\chi^{\xi^{(1)}(\xi^{-1})^{(1)}},\\
\Ind_{PGO_2^-(\Fq)}^{PGL_2(\Fq)}(1)&=\chi^{1^{(2)}}
+\left\{\begin{array}{cl}\chi^{\eta^{(1^2)}},&\text{if }d_\eta=1\\0,
&\text{if }d_\eta=-1\end{array}\right\}\\
&\qquad\qquad+\sum_{\substack{\xi\in L^\sigma\setminus\{1,\eta\}\\d_\xi=1}}
\chi^{\xi^{(1)}(\xi^{-1})^{(1)}}
+\sum_{\substack{\xi\in L^\tsigma\setminus\{1,\eta\}\\\tilde{d}_\xi=-1}}
\chi^{\xi^{(1)}(\xi^{-1})^{(1)}}.
\end{split}
\eeq
In the sums, only one representative of each pair $\{\xi,\xi^{-1}\}$ should be taken.
Of course, these
decompositions are known and easy to calculate directly.
(As Theorem \ref{firstpgothm} makes clear, the fact that they are multiplicity-free
is specific to the case $n=2$.)
\eex

As mentioned in the introduction, we will prove Theorems
\ref{firstpgspthm} and \ref{firstpgothm} in the equivalent forms involving
the basic characters $B_{\unnu}$. Applying \eqref{glneqn}, these are:
\bth \label{pgspthm}
For any $\unnu\in\widehat{\bbP}_n^{\sigma}$,
\[ \langle B_{\unnu},
\Ind_{PGSp_n(\Fq)}^{PGL_n(\Fq)}(1)\rangle =
\left\{\begin{array}{cl}{\displaystyle
\prod_{\xi\in\subsigL}\sum_{\substack{\rho_\xi\vdash|\nu_\xi|\\
\rho_\xi\textup{ even}}}
\chi_{\nu_\xi}^{\rho_\xi}},&\begin{array}{c}\text{if all $|\nu_\xi|$ are even and}\\
\prod_{\xi\in\langle\sigma\rangle\setminus L} 
N(\xi)^{\frac{|\nu_\xi|}{2}}=1,\end{array}\\
0, &\text{ otherwise.}
\end{array}\right. \]
\eth
\bth \label{pgothm}
For any $\unnu\in\widehat{\bbP}_n^{\sigma}$ and $\epsilon\in\{\pm 1\}$,
\bes
\begin{split}
\langle B_{\unnu},
\Ind_{PGO_n^\epsilon(\Fq)}^{PGL_n(\Fq)}(1)\rangle &=
\frac{1}{4}\prod_{\substack{\xi\in\langle\sigma\rangle\setminus L\\
d_\xi=1}}(-1)^{|\nu_\xi|}\sum_{\rho_\xi\vdash|\nu_\xi|}
(\prod_i (m_i(\rho_\xi)+1))\chi_{\nu_\xi}^{\rho_\xi}\\
&\qquad\thickspace\thickspace\times
\prod_{\substack{\xi\in\langle\sigma\rangle\setminus L\\ d_\xi=-1}}
\sum_{\substack{\rho_\xi\vdash|\nu_\xi|\\\rho_\xi' \textup{ even}}}
\chi_{\nu_\xi}^{\rho_\xi}\\
&\negthickspace\negthickspace\negthickspace\negthickspace
\negthickspace\negthickspace\negthickspace\negthickspace
\negthickspace\negthickspace\negthickspace\negthickspace
+\left\{\begin{array}{cl}{\displaystyle
\frac{1}{2}\epsilon\prod_{\xi\in\subsigL}
\sum_{\substack{\rho_\xi\vdash|\nu_\xi|\\\rho_\xi'\textup{ even}}}
\chi_{\nu_\xi}^{\rho_\xi}},&\begin{array}{c}\text{if all $|\nu_\xi|$ are even and}\\ 
\prod_{\xi\in\langle\sigma\rangle\setminus L} 
N(\xi)^{\frac{|\nu_\xi|}{2}}=1,\end{array}\\
0, &\text{ otherwise}
\end{array}\right.\\
&\negthickspace\negthickspace\negthickspace\negthickspace
\negthickspace\negthickspace\negthickspace\negthickspace
\negthickspace\negthickspace\negthickspace\negthickspace
+\frac{\Phi(\unnu)}{4}
\prod_{\substack{\xi\in\langle\sigma\rangle\setminus L\\d_\xi=1\\2\nmid m_\xi}}
\sum_{\substack{\rho_\xi\vdash|\nu_\xi|\\2\mid m_{2i+1}(\rho_\xi), \forall i}}
(-1)^{\frac{1}{2}|\rho_\xi|+\ell(\rho_\xi)_{2\textup{ mod }4}}
(\prod_i (m_{2i}(\rho_\xi)+1))\chi_{\nu_\xi}^{\rho_\xi}\\
&\times
\prod_{\substack{\xi\in\langle\sigma\rangle\setminus L\\d_\xi=1\\2\mid m_\xi}}
(-1)^{|\nu_\xi|+\frac{1}{2}m_\xi|\nu_\xi|}\sum_{\rho_\xi\vdash|\nu_\xi|}
(\prod_i (m_i(\rho_\xi)+1))\chi_{\nu_\xi}^{\rho_\xi}\\
&\times
\prod_{\substack{\xi\in\langle\sigma\rangle\setminus L\\ d_\xi=-1}}
(-1)^{\frac{1}{2}m_\xi|\nu_\xi|}
\sum_{\substack{\rho_\xi\vdash|\nu_\xi|\\\rho_\xi' \textup{ even}}}
\chi_{\nu_\xi}^{\rho_\xi}.
\end{split}
\ees
\eth
\noindent
Here the sign $(-1)^{n+\sum_{\xi\in\langle\sigma\rangle\setminus L} |\nu_\xi|}$
of \eqref{glneqn} has been rewritten as $\prod_{\xi\in\subsigL}(-1)^{|\nu_\xi|}$
and then distributed among the factors as appropriate. Also, in the last term of
Theorem \ref{pgothm}, the sign $(-1)^{\frac{n}{2}}$ of Theorem \ref{firstpgothm}
has been rewritten as $\prod_{\xi\in\subsigL}(-1)^{\frac{1}{2}m_\xi|\nu_\xi|}$
and then distributed among the factors (each of which can only be nonzero
when $m_\xi|\nu_\xi|$ is even, so this and $\Phi(\unnu)$ make sense).
Theorems \ref{pgspthm} and \ref{pgothm} will be deduced from Lusztig's formula
in the next two sections.
\section{Proof of Theorem \ref{pgspthm}}
The proof of Theorem \ref{pgspthm} is mostly identical to that of
\cite[(2.1.1)]{mysymm}, but we will go through all the steps
again as a warm-up for the more complicated arguments in the next section.
We fix a non-degenerate \textbf{skew-symmetric} bilinear form 
$B$ on $V$, compatible with
$F$ in the sense that $B(F(v),F(v')) = B(v,v')^q$. 
Let $\theta:G\to G$ be the involution defined by
\[ B(\theta(g)v,v')=B(v,g^{-1}v'),\ \forall g\in G,
v,v'\in V, \]
and let $K=G^\theta=Sp(V,B)$.
We have $F\theta=\theta F$, so $K$ is $F$-stable, and 
$K^F$ is what we have been calling $Sp_n(\Fq)$.

Now for $z\in Z$, $\theta(z)=z^{-1}$, so $K\cap Z=\{\pm 1\}$.
If we write $\theta$ also for the induced involution of $\bG$, it is
clear that $\bG^\theta=\bK=KZ/Z$, since
\[ KZ=\{g\in G\,|\,B(gv,gv')=\alpha B(v,v'),\
\forall v,v'\in V,\,\exists\alpha\in k^\times\}, \] 
the group called either $GSp_n$ or $CSp_n$. However, as mentioned in the
introduction,
\[ K^F Z^F=\{g\in G^F\,|\,B(gv,gv')=\alpha B(v,v'),\
\forall v,v'\in V,\,\exists\alpha\in (\Fq^\times)^2\} \]
is a subgroup of index $2$ in
\[ (KZ)^F =\{g\in G^F\,|\,B(gv,gv')=\alpha B(v,v'),\
\forall v,v'\in V,\,\exists\alpha\in \Fq^\times\}, \]
so the image of $K^F\to\bK^F$ is of index $2$ in $\bK^F$. The group $\bK^F$
is what we have been calling $PGSp_n(\Fq)$.

Now let $\unnu\in\widehat{\bbP}_n^{\sigma}$, and fix an $F$-stable
maximal torus $\bT$ of $\bG$ and a character $\lambda\in\widehat{\bT^F}$
such that the pair $(\bT,\lambda)$ is in the $\bG^F$-orbit corresponding
to $\unnu$ under the bijection described in the previous section.
Let $T$ be the preimage of $\bT$ in $G$, a maximal torus of $G$.
We want to deduce Theorem \ref{pgspthm} from Lusztig's formula 
\cite[Theorem 3.3]{symmfinite} for the inner product
$\langle B_{\unnu},\Ind_{\bK^F}^{\bG^F}(1)\rangle$. A key ingredient
in the formula is the set
\[ \Theta_{\bT}:=\{\bar{f}\in\bG\,|\,\theta(\bar{f}^{-1}\bT\bar{f})
=\bar{f}^{-1}\bT\bar{f}\}. \]
Since $T\supset Z$, it is obvious that 
$\Theta_{\bT}$ is the image of the analogous
set 
\[ \Theta_T:=\{f\in G\,|\,\theta(f^{-1}Tf)=f^{-1}Tf\}. \] 
For $f\in\Theta_T$, $f^{-1}Tf$ is a $\theta$-stable
maximal torus of $G$, or in other words the stabilizer of a decomposition
$V=\bigoplus_{i=1}^n L_i$ into lines with the property that
$L_j^\perp=\bigoplus_{i\neq w(j)} L_i$ for some
involution $w\in S_n$ with no fixed points. (Here $L_j^\perp$ is the
subspace perpendicular to $L_j$ under $B$.)
If $t\in f^{-1}Tf$ has eigenvalue $\alpha_i$ on $L_i$ for all $i$, then
$\theta(t)$ has eigenvalue $\alpha_{w(i)}^{-1}$ on $L_i$ for all $i$;
thus $t$ lies in $f^{-1}Tf\cap K$ if and only if
$\alpha_i\alpha_{w(i)}=1$ for all $i$, and
the image $\bar{t}$ of $t$ lies in $\bar{f}^{-1}\bT\bar{f}\cap\bK$
if and only if $\alpha_i\alpha_{w(i)}=\beta$ for all $i$, for some
$\beta\in k^\times$.
From this the following result is obvious.
\bpr \label{connprop}
For all $\bar{f}\in\Theta_{\bT}$, $\bar{f}^{-1}\bT\bar{f}\cap\bK$ is
connected and $Z_{\bG}(\bar{f}^{-1}\bT\bar{f}\cap\bK)=\bar{f}^{-1}\bT\bar{f}$.
\epr
\noindent
Under these circumstances, Lusztig's formula takes a
simpler form. Define 
\[ \Theta_{\bT,\lambda}^F:=\{\bar{f}\in
\Theta_{\bT}^F\,|\, 
\lambda\vert_{(\bT\cap\bar{f}^{-1}\bK\bar{f})^{F}}=1\}. \]
This is clearly a union of $\bT^F$--\,$\bK^F$ double cosets, and
\cite[Theorem 3.3]{symmfinite} implies that
\beq \label{simplelusztigeqn}
\langle B_{\unnu},\Ind_{\bK^F}^{\bG^F}(1)\rangle
=|\bT^F\leftdiv\Theta_{\bT,\lambda}^F/\bK^F|.
\eeq
We now aim to express the right-hand side in terms of $\unnu$.

Let $W(T)=N_G(T)/T$ be the Weyl group of $T$ (isomorphic to the symmetric
group $S_n$), and let $W(T)_{\ffi}$ be the subset corresponding to the
fixed-point-free involutions in $S_n$. 
For any $f\in\Theta_T$, we define $w_{f}\in W(T)_{\ffi}$
by requiring that $f^{-1}w_f f\in W(f^{-1}Tf)$ corresponds to the
above involution $w$ of the fixed lines of $f^{-1}Tf$. It is clear that
$w_f$ depends only on $\bar{f}\in\Theta_{\bT}$.
\bpr \label{geombij}
The map $\bar{f}\mapsto w_f$ induces a bijection
$\bT\leftdiv\Theta_{\bT}/\bK\isomto W(T)_{\ffi}$.
\epr
\bpf
This follows trivially from the corresponding statement for $G$, which is
\cite[Proposition 2.0.1]{mysymm}.
\epf
\bpr \label{fbij}
The map $\bar{f}\mapsto w_f$ induces a bijection
$\bT^F\leftdiv\Theta_{\bT}^F/\bK^F\isomto W(T)_{\ffi}^F$.
\epr
\bpf
Clearly the map in Proposition \ref{geombij} is $F$-stable, so we get
a bijection 
\[ (\bT\leftdiv\Theta_{\bT}/\bK)^F\isomto W(T)_{\ffi}^F. \]
Now regard $\bT\leftdiv\Theta_{\bT}/\bK$ as the set of orbits
of the group $\bT\times\bK$ on $\Theta_{\bT}$. Since
$\bT\times\bK$ is connected, Lang's Theorem implies that
$\bT^F\leftdiv\Theta_{\bT}^F/\bK^F\to (\bT\leftdiv\Theta_{\bT}/\bK)^F$
is surjective. To prove injectivity, we need only check that the stabilizers
in $\bT\times\bK$ of points in $\Theta_{\bT}$ are connected, which is
exactly what the first part of Proposition \ref{connprop} says.
(Note that in the proof of the corresponding result \cite[Lemma 2.1.2]{mysymm},
the words `injective' and `surjective' need to be swapped!)
\epf

The next step is to identify the subset of $W(T)_{\ffi}^F$ which
corresponds under this bijection to the set we are trying to enumerate,
$\bT^F\leftdiv\Theta_{\bT,\lambda}^F/\bK^F$. Define 
\[ W(T)_{\lambda,\ffi}^F:=\{w\in W(T)_{\ffi}^F\,|\,\lambda({}^{w}t)=\lambda(t),
\,\forall t\in T^F\}. \]
(Here we are abusing notation by writing $\lambda$ also for the pull-back
of $\lambda$ to $T^F$.) Now \cite[Lemma 2.1.3]{mysymm} states that for
$f\in\Theta_T^F$, $f\in\Theta_{T,\lambda}^F\Longleftrightarrow
w_f\in W(T)_{\lambda,\ffi}^F$. This follows immediately from the
fact that
\beq 
(T\cap fKf^{-1})^F=\{t\in T^F\,|\,{}^{w_f}t=t^{-1}\}=
\{({}^{w_f}t)t^{-1}\,|\,t\in T^F\}.
\eeq
However, for $\bar{f}\in\Theta_{\bT}^F$ the analogous result is not quite true,
since
\[ \bT_{\bar{f}}^F:=\{({}^{w_f}\bar{t})\bar{t}^{-1}\,|\,\bar{t}\in\bT^F\} \] 
is a subgroup
of index $2$ in 
\[ (\bT\cap \bar{f}\bK\bar{f}^{-1})^F=\{\bar{t}\in\bT^F\,|\,
{}^{w_f}\bar{t}=\bar{t}^{-1}\}. \] 
To see this, recall that the fixed lines
of $T$ can be labelled $L_{(\xi,j,i)}$ as in \eqref{toruseqn}, so
$W(T)$ can be viewed as the group of permutations of the triples $(\xi,j,i)$.
If $t\in T^F$ has eigenvalue $\alpha_{(\xi,j,i)}$ on $L_{(\xi,j,i)}$, then
\beq 
{}^{w_f}\bar{t}=\bar{t}^{-1}\Longleftrightarrow \alpha_{(\xi,j,i)}
\alpha_{w_f(\xi,j,i)}=\beta,\,\forall (\xi,j,i),\text{ for some
$\beta\in\Fq^\times$,}
\eeq
whereas
\beq
\bar{t}\in \bT_{\bar{f}}^F\Longleftrightarrow \alpha_{(\xi,j,i)}
\alpha_{w_f(\xi,j,i)}=\beta,\,\forall (\xi,j,i),\text{ for some
$\beta\in(\Fq^\times)^2$.}
\eeq
Consequently, the correct analogue of \cite[Lemma 2.1.3]{mysymm} is:
\blm \label{unsatislemma}
For $\bar{f}\in\Theta_{\bT}^F$,
$\bar{f}\in\Theta_{\bT,\lambda}^F$ if and only if
$w_f\in W(T)_{\lambda,\ffi}^F$ \textbf{and} $\lambda(t_0)=1$ for some
(hence any) $t_0\in T^F$ whose eigenvalues $\{\alpha_{(\xi,j,i)}\}$
have the property that $\alpha_{(\xi,j,i)}
\alpha_{w_f(\xi,j,i)}=\beta,\,\forall (\xi,j,i)$, for some
$\beta\in\Fq^\times\setminus(\Fq^\times)^2$.
\elm 

As in \cite[\S1.4]{mysymm}, the set $W(T)_{\lambda,\ffi}^F$ can be described
using \eqref{toruseqn}--\eqref{lambdaeqn}. 
Let $\Lambda$ be the set of pairs $(\xi,j)$ where $\xi\in\sigL$ and 
$1\leq j\leq\ell(\nu_\xi)$. (In \cite{mysymm}, this and the related sets
were denoted $\Lambda(\unnu)$, etc.) Any
$w\in W(T)_{\lambda,\ffi}^F$ defines a partition $\Lambda=\Lambda_w^2\coprod
\Lambda_w^3$ such that:
\begin{itemize}
\item for all $(\xi,j)\in\Lambda_w^2$, $(\nu_\xi)_j$ is even and
\[ w(\xi,j,i)=(\xi,j,i+\frac{1}{2}m_\xi(\nu_\xi)_j)\text{ for all }
i\in\Z/m_\xi(\nu_\xi)_j\Z; \]
\item for all $(\xi,j)\in\Lambda_w^3$, 
\[ w(\xi,j,i)=(\xi,\tilde{j},i+i(w,\xi,j))\text{ for all }
i\in\Z/m_\xi(\nu_\xi)_j\Z, \]
where $(\xi,\tilde{j})\in\Lambda_w^3$, $\tilde{j}\neq j$,
$(\nu_\xi)_{\tilde{j}}=(\nu_\xi)_j$, and
$i(w,\xi,j)$ is some fixed `shift' in $m_\xi\Z/m_\xi(\nu_\xi)_j\Z$.
\end{itemize}
(The set $\Lambda_w^1$ will appear in the next section; it is empty here
because we are dealing with fixed-point-free involutions.)
Clearly $(\xi,j)\mapsto(\xi,\tilde{j})$ is a fixed-point-free involution
of $\Lambda_w^3$; we will write $\widetilde{\Lambda_w^3}$ for an arbitrarily chosen
set of representatives of its orbits.
Note that the existence of $w\in W(T)_{\lambda,\ffi}^F$ forces each 
$|\nu_\xi|$ to be even.

As explained in \cite[\S1.4, \S2.1]{mysymm}, the above description gives
a bijection between $W(T)_{\lambda,\ffi}^F$ and 
$\prod_{\xi\in\subsigL} Z_{\ffi}^{\nu_\xi}$, where
$Z_{\ffi}^{\nu}$ denotes the set of fixed-point-free involutions in
$S_{|\nu|}$ which commute with a given element $w_\nu$ of cycle-type $\nu$;
essentially this bijection is defined by dividing the shifts
$\frac{1}{2}m_\xi(\nu_\xi)_j$ and $i(w,\xi,j)$ by $m_\xi$.

Now suppose $t_0$, $\alpha_{(\xi,j,i)}$, $\beta$ are as in Lemma
\ref{unsatislemma}. Since $t_0\in T^F$, we have $\alpha_{(\xi,j,i+1)}
=\alpha_{(\xi,j,i)}^q$ for all $(\xi,j,i)$. Clearly the condition
$\alpha_{(\xi,j,i)}\alpha_{w_f(\xi,j,i)}=\beta$ is equivalent to the
conjunction of the following conditions:
\begin{itemize}
\item for all $(\xi,j)\in\Lambda_{w_f}^2$, 
$\alpha_{(\xi,j,1)}^{1+q^{\frac{1}{2}m_\xi(\nu_\xi)_j}}=\beta$;
\item for all $(\xi,j)\in\Lambda_{w_f}^3$,
$\alpha_{(\xi,j,1)}\alpha_{(\xi,\tilde{j},1)}^{q^{i(w_f,\xi,j)}}=\beta$.
\end{itemize}
Hence
\bes
\begin{split}
\lambda(t_0)&=\prod_{\xi\in\subsigL}\prod_{j=1}^{\ell(\nu_\xi)}
\lng\alpha_{(\xi,j,1)},\xi\rng_{m_\xi(\nu_\xi)_j}\\
&=\prod_{(\xi,j)\in\Lambda_{w_f}^2}
\lng\alpha_{(\xi,j,1)}^{1+q^{m_\xi}+q^{2m_\xi}+\cdots+q^{((\nu_\xi)_j-1)m_\xi}}
,\xi\rng_{m_\xi}\\
&\qquad\times
\prod_{(\xi,j)\in\widetilde{\Lambda_{w_f}^3}}
\lng\alpha_{(\xi,j,1)}^{1+q^{m_\xi}+q^{2m_\xi}+\cdots+q^{((\nu_\xi)_j-1)m_\xi}}
\alpha_{(\xi,\tilde{j},1)}
^{1+q^{m_\xi}+q^{2m_\xi}+\cdots+q^{((\nu_\xi)_{\tilde{j}}-1)m_\xi}}
,\xi\rng_{m_\xi}\\
&=\prod_{(\xi,j)\in\Lambda_{w_f}^2}
\lng\beta^{\frac{1}{2}(\nu_\xi)_j},\xi\rng_{m_\xi}
\prod_{(\xi,j)\in\widetilde{\Lambda_{w_f}^3}}
\lng\beta^{(\nu_\xi)_j},\xi\rng_{m_\xi}\\
&=\prod_{(\xi,j)\in\Lambda_{w_f}^2}
\lng\beta,N(\xi)\rng^{\frac{1}{2}(\nu_\xi)_j}
\prod_{(\xi,j)\in\widetilde{\Lambda_{w_f}^3}}
\lng\beta,N(\xi)\rng^{(\nu_\xi)_j}\\
&=\lng\beta,\prod_{\xi\in\subsigL}N(\xi)^{\frac{1}{2}|\nu_\xi|}\rng.
\end{split}
\ees
Note that $\prod_{\xi\in\subsigL}N(\xi)^{\frac{1}{2}|\nu_\xi|}$ 
is a square root
of $\Pi(\unnu)=1$.
So Lemma \ref{unsatislemma} can be restated as follows:
\bpr \label{spbijprop}
If all $|\nu_\xi|$ are even and
$\prod_{\xi\in\subsigL}N(\xi)^{\frac{1}{2}|\nu_\xi|}=1$,
$\bT^F\leftdiv\Theta_{\bT,\lambda}^F/\bK^F$ is in bijection with
$W(T)_{\lambda,\ffi}^F$. Otherwise,
$\Theta_{\bT,\lambda}^F$ is empty. 
\epr
\noindent
Combining this with Lusztig's formula \eqref{simplelusztigeqn}, we deduce
\bpr
For any $\unnu\in\widehat{\bbP}_n^{\sigma}$,
\[ \langle B_{\unnu},\Ind_{\bK^F}^{\bG^F}(1)\rangle
=\left\{
\begin{array}{cl}
\prod_{\xi\in\subsigL}|Z_{\ffi}^{\nu_\xi}|,&\begin{array}{c}\text{if all 
$|\nu_\xi|$ are even and}\\
\prod_{\xi\in\subsigL}N(\xi)^{\frac{1}{2}|\nu_\xi|}=1,\end{array}\\
0,&\text{ otherwise.}
\end{array}
\right. \]
\epr
\noindent
This implies Theorem \ref{pgspthm} by the same combinatorial fact used in
\cite[\S2.1]{mysymm}: namely,
\beq \label{macdonaldeqn}
|Z^{\nu}_{\ffi}|
=\sum_{\substack{\rho\vdash|\nu|\\\rho\text{ even}}}
\chi_\nu^{\rho}.
\eeq
A reference for this fact is \cite[VII(2.4)]{macdonald}.
\section{Proof of Theorem \ref{pgothm}}
In this section we fix a non-degenerate \textbf{symmetric} bilinear form 
$B$ on $V$, compatible with
$F$ in the sense that $B(F(v),F(v')) = B(v,v')^q$. 
Let $\theta:G\to G$ be the involution defined by
\[ B(\theta(g)v,v')=B(v,g^{-1}v'),\ \forall g\in G,
v,v'\in V, \]
and let $K=G^\theta=O(V,B)$. Note that $K$
has two components, with $K^\circ=SO(V,B)$;
since $K\cap Z=\{\pm 1\}\subset K^\circ$, $\bG^\theta=\bK=KZ/Z$ has two
components also. This group $\bK$ is what is known as $PGO_n$.

We have $F\theta=\theta F$, so $K$ is $F$-stable.
The Witt index of $B$ on $V^F$ is either $n/2$
or $n/2-1$, and accordingly $\bK^F\cong PGO_n^\epsilon(\Fq)$ where
$\epsilon=+1$ or $-1$. As in the previous section, the image of $K^F\to\bK^F$
is an index-$2$ subgroup of $\bK^F$ (\textbf{not} equal to
the index-$2$ subgroup $(\bK^\circ)^F$).

Fix $\unnu\in\widehat{\bbP}_n^{\sigma}$, and define $\bT,\lambda,T,
\Theta_{\bT},\Theta_T,W(T)\cong S_n$ as in the previous section.
Let $W(T)_\inv$ denote the set of involutions in $W(T)$ (including $1$).
For any $f\in\Theta_T$, $f^{-1}Tf$ is the stabilizer of a decomposition
$V=\bigoplus_{i=1}^n L_i$ into lines with the property that
$L_j^\perp=\bigoplus_{i\neq w(j)} L_i$ for some
involution $w\in S_n$, and we define $w_f\in W(T)_\inv$ by requiring that
$f^{-1}w_f f\in W(f^{-1}Tf)$ corresponds to $w$. It is clear that
$w_f$ depends only on $\bar{f}\in\Theta_{\bT}$.
\bpr
The map $\bar{f}\mapsto w_f$ induces a bijection
$\bT\leftdiv\Theta_{\bT}/\bK\isomto W(T)_{\inv}$.
\epr
\bpf
This follows trivially from the corresponding statement for $G$, which is
\cite[Proposition 4.0.2]{mysymm}.
\epf
However, the analogue of Proposition \ref{fbij} is false, since neither
$\bK$ nor $\bT\cap\bar{f}\bK\bar{f}^{-1}$ for general 
$\bar{f}\in\Theta_{\bT}^F$ is connected. Insead, we have the following result,
where $\epsilon_\unnu$ denotes the sign of the permutation by which $F$ acts 
on the fixed lines of $T$:
\bpr \label{noninjmapprop}
The map $\bar{f}\mapsto w_f$ induces a map
$\bT^F\leftdiv\Theta_{\bT}^F/\bK^F\to W(T)_{\inv}^F$.
If $\epsilon_\unnu=\epsilon$, this map is surjective. If
$\epsilon_\unnu=-\epsilon$, the image is
$W(T)_{\inv}^F\setminus W(T)_{\ffi}^F$.
\epr
\bpf
This follows trivially from the corresponding statement for $G^F$, which is
\cite[Lemma 4.2.2]{mysymm}.
\epf

To state Lusztig's formula in this case, we need to introduce, for any
$\bar{f}\in\Theta_{\bT}^F$, the function 
$\epsilon_{\bT,\bar{f}}:(\bT\cap\bar{f}\bK\bar{f}^{-1})^F\to\{\pm 1\}$
defined by
\[ \epsilon_{\bT,\bar{f}}(\bar{t})=(-1)^{
\text{$\Fq$-rank}(Z_{\bG}((\bT\cap \bar{f}\bK\bar{f}^{-1})^{\circ}))
+\text{$\Fq$-rank}(Z_{\bG}^{\circ}(\bar{t})\cap Z_{\bG}((\bT\cap 
\bar{f}\bK\bar{f}^{-1})^{\circ}))}. \]
It follows from \cite[Proposition 2.3]{symmfinite} that
$\epsilon_{\bT,\bar{f}}$ is a group homomorphism which factors through
$(\bT\cap \bar{f}\bK\bar{f}^{-1})^{F}/
((\bT\cap \bar{f}\bK\bar{f}^{-1})^{\circ})^{F}$. We can then define
\[ \Theta_{\bT,\lambda}^{F}:=\{\bar{f}\in\Theta_{\bT}^{F}\,|\,
\lambda|_{(\bT\cap \bar{f}\bK\bar{f}^{-1})^{F}}=\epsilon_{\bT,\bar{f}}\}. \]
This is clearly a union of $\bT^{F}$--\,$\bK^{F}$ double cosets,
and \cite[Theorem 3.3]{symmfinite} says that
\beq \label{lusztigeqn}
\langle B_\unnu, \Ind_{\bK^{F}}^{\bG^{F}}(1)\rangle
=\sum_{\bar{f}\in \bT^{F}\setminus \Theta_{\bT,\lambda}^{F}/\bK^{F}}
(-1)^{\textup{$\Fq$-rank}(\bT)+\textup{$\Fq$-rank}
(Z_{\bG}((\bT\cap \bar{f}\bK\bar{f}^{-1})^{\circ}))}.
\eeq
Our immediate aim is to use the map $\bar{f}\mapsto w_f$ to 
turn the right-hand side
into a sum over a suitable subset of $W(T)_\inv^F$.

As in the previous section, we view $W(T)$ as the group of permutations
of the triples $(\xi,j,i)$ in \eqref{toruseqn}. If $\bar{f}\in\Theta_{\bT}$,
and $t\in T$ has eigenvalue $\alpha_{(\xi,j,i)}$ on $L_{(\xi,j,i)}$, then
\beq \label{descripeqn} 
\bar{t}\in \bT\cap \bar{f}\bK\bar{f}^{-1}
\Longleftrightarrow
\alpha_{(\xi,j,i)}\alpha_{w_f(\xi,j,i)}=\beta,\,\forall (\xi,j,i),
\text{ for some $\beta\in k^\times$.}
\eeq
Unlike in the previous section, $w_f$ is now allowed to have fixed points;
the eigenvalues $\alpha_{(\xi,j,i)}$ where $w_f(\xi,j,i)=(\xi,j,i)$ must all be
square roots of $\beta$. It is clear that the additional condition required
in order that $\bar{t}\in (\bT\cap \bar{f}\bK\bar{f}^{-1})^\circ$ is that
all these square roots are equal.

This allows us to prove a partial analogue of Lemma \ref{unsatislemma}. Define
\[ W(T)_{\lambda,\inv}^F:=\{w\in W(T)_\inv^F\,|\,\lambda({}^w t)=\lambda(t),
\,\forall t\in T^F\}. \]
(Once again we are blurring the distinction between $\lambda$ 
and its pull-back to $T^F$.)
\blm \label{partiallemma} 
If $\bar{f}\in\Theta_{\bT,\lambda}^F$, then
$w_f\in W(T)_{\lambda,\inv}^F$.
\elm
\bpf
From our description of $(\bT\cap \bar{f}\bK\bar{f}^{-1})^\circ$, it is clear that
for all $t\in T$ the element
$({}^{w_f}\bar{t})\bar{t}^{-1}$ lies in 
$(\bT\cap \bar{f}\bK\bar{f}^{-1})^\circ$. Since $\epsilon_{\bT,\bar{f}}$
is trivial on $((\bT\cap \bar{f}\bK\bar{f}^{-1})^{\circ})^{F}$, the assumption
$\bar{f}\in\Theta_{\bT,\lambda}^F$ implies that
$\lambda(({}^{w_f}\bar{t})\bar{t}^{-1})=1$ for all $t\in T^F$.
Thus $w_f\in W(T)_{\lambda,\inv}^F$.
\epf

The set $W(T)_{\lambda,\inv}^F$ contains the set
$W(T)_{\lambda,\ffi}^F$ used in the previous section and can be described 
similarly. Namely, any $w\in W(T)_{\lambda,\inv}^F$ defines a partition
$\Lambda=\Lambda_w^1\coprod\Lambda_w^2\coprod\Lambda_w^3$ such that
for $(\xi,j)\in\Lambda_w^1$,
$w(\xi,j,i)=(\xi,j,i)$ for all $i\in\Z/m_\xi(\nu_\xi)_j\Z$, and $\Lambda_w^2$
and $\Lambda_w^3$ are as before. Define $\widetilde{\Lambda_w^3}$ as before also,
and let $\ell_w^i=|\Lambda_w^i|$.  Note that since $n$ is even,
\beq \label{eveneqn}
\sum_{(\xi,j)\in\Lambda_w^1} m_\xi(\nu_\xi)_j\text{ is even.}
\eeq

As explained in \cite[\S1.4]{mysymm}, this description gives
a bijection between $W(T)_{\lambda,\inv}^F$ and 
$\prod_{\xi\in\subsigL} Z_{\inv}^{\nu_\xi}$, where
$Z_{\inv}^{\nu}$ denotes the set of fixed-point-free involutions in
$S_{|\nu|}$ which commute with a given element $w_\nu$ of cycle-type $\nu$;
this bijection is defined by dividing all shifts by $m_\xi$.

Now we define a subset $X$ of
$W(T)_{\lambda,\inv}^F$ by the rule that $w\in X$ if and only if
$(\nu_\xi)_j$ is even for all $(\xi,j)\in\Lambda_w^1$ with $d_\xi=-1$.
(In \cite[\S4]{mysymm} this was called $X_\inv^{\unnu}$.)
Also let $Y$ be the subset of
$W(T)_{\lambda,\inv}^F$ defined by the requirement that $m_\xi(\nu_\xi)_j$
is even for all $(\xi,j)\in\Lambda_w^1$. Trivially
$W(T)_{\lambda,\ffi}^F\subseteq X\cap Y$.
We can measure the non-injectivity of the map in Proposition \ref{noninjmapprop}:
\bpr \label{fibreprop}
Let $w\in W(T)_{\lambda,\inv}^F$. The number of double cosets $\bT^F\bar{f}\,\bK^F$
in $\bT^F\leftdiv\Theta_{\bT}^F/\bK^F$ such that $w_f=w$ equals
\[ \left\{\begin{array}{cl}
0,&\text{ if $w\in W(T)_{\lambda,\ffi}^F$ and $\epsilon_{\unnu}=-\epsilon$,}\\
1,&\text{ if $w\in W(T)_{\lambda,\ffi}^F$ and $\epsilon_{\unnu}=\epsilon$,}\\
2^{\ell_w^1-1},&\text{ if $w\in Y\setminus W(T)_{\lambda,\ffi}^F$,}\\
2^{\ell_w^1-2},&\text{ if $w\in W(T)_{\lambda,\inv}^F\setminus Y$.}
\end{array}\right. \]
\epr
\bpf
Let $\bar{f}\in\Theta_{\bT}^F$ be such that $w_f\in W(T)_{\lambda,\inv}^F$.
The number of $\bT^F$--\,$\bK^F$ double cosets in $(\bT\bar{f}\bK)^F$ can be
calculated by the same method as in \cite[Lemma 4.1.3]{mysymm}: 
it is the index in the 
group
\[ H:=\left\{(\epsilon_{(\xi,j,i)})\in
\negthickspace\negthickspace
\prod_{(\xi,j)\in\Lambda_{w_f}^{1}}\negthickspace
(\pm 1)^{m_\xi(\nu_\xi)_j}\,\left|\,
\prod_{(\xi,j,i)}\epsilon_{(\xi,j,i)}=1\right.\right\} \]
of the subgroup $H'$ defined by the further requirement that either
$\prod_i \epsilon_{(\xi,j,i)} =1$ for all $(\xi,j)\in\Lambda_{w_f}^{1}$
or $\prod_i \epsilon_{(\xi,j,i)} =(-1)^{m_\xi(\nu_\xi)_j}$ 
for all $(\xi,j)\in\Lambda_{w_f}^{1}$. If $w_f\in W(T)_{\lambda,\ffi}^F$, then
$\Lambda_{w_f}^{1}$ is empty and $H$ is trivial, so the index is $1$. Otherwise,
\beq 
|H|=2^{(\sum_{(\xi,j)\in\Lambda_{w_f}^{1}} m_\xi(\nu_\xi)_j)-1}. 
\eeq
If $w_f\in Y\setminus W(T)_{\lambda,\ffi}^F$, then the two possibilities 
in the definition of $H'$ are the same, so $|H'|=2^{\sum(m_\xi(\nu_\xi)_j-1)}$
and the index is $2^{\ell_{w_f}^1-1}$. If $w_f\notin Y$, then $H'$ is twice
as large (bearing in mind \eqref{eveneqn}), so the index is
$2^{\ell_{w_f}^1-2}$. Combining this with Proposition \ref{noninjmapprop},
we have the result.
\epf

To complete the interpretation of \eqref{lusztigeqn}, we define, for each
$w\in Y$, the sign
\[ \Phi(w):= \prod_{(\xi,j)\in\Lambda_w^1}
(-1)^{\frac{1}{2}m_\xi(\nu_\xi)_j}
\prod_{(\xi,j)\in\Lambda_w^2}
d_\xi^{\frac{1}{2}m_\xi(\nu_\xi)_j}
\prod_{(\xi,j)\in\widetilde{\Lambda_w^3}}
d_\xi^{m_\xi(\nu_\xi)_j}. \]
\bpr \label{mainprop}
Let $\bar{f}\in\Theta_{\bT}^F$ be such that $w_f\in W(T)_{\lambda,\inv}^F$.
\ben
\item $\textup{$\Fq$-rank}(\bT)+\textup{$\Fq$-rank}
(Z_{\bG}((\bT\cap \bar{f}\bK\bar{f}^{-1})^{\circ}))\equiv\ell_{w_f}^1$ 
\textup{(mod $2$)}.
\item $\bar{f}\in\Theta_{\bT,\lambda}^F$ if and only if \textbf{either}
$w_f\in X\setminus Y$ \textbf{or} $w_f\in X\cap Y$ and $\Phi(w_f)=\Phi(\unnu)$.
\een
\epr
\bpf
The $\Fq$-rank of $\bT$ is $\ell-1$, where $\ell=|\Lambda|=\sum_{\xi\in\subsigL}
\ell(\nu_\xi)$ is the total number of pairs $(\xi,j)$. Using the description
of $(\bT\cap\bar{f}\bK\bar{f}^{-1})^\circ$ given after \eqref{descripeqn}
we see that
\beq 
Z_{\bG}((\bT\cap\bar{f}\bK\bar{f}^{-1})^\circ)=
(GL(\negthickspace\negthickspace\bigoplus_{\substack{(\xi,j,i)\\
(\xi,j)\in\Lambda_{w_f}^{1}}}\negthickspace L_{(\xi,j,i)})\
\times\negthickspace
\prod_{\substack{(\xi,j,i)\\(\xi,j)\notin\Lambda_{w_f}^{1}}}\negthickspace
GL(L_{(\xi,j,i)}))/Z,
\eeq
which has $\Fq$-rank
\[ \sum_{(\xi,j)\in\Lambda_{w_f}^{1}} m_\xi(\nu_\xi)_j +\ell-\ell_{w_f}^1-1, \]
and (1) follows using \eqref{eveneqn}. To prove (2), recall that by definition
$\bar{f}\in\Theta_{\bT,\lambda}^F$ iff $\epsilon_{\bT,\bar{f}}(\bar{t})=
\lambda(\bar{t})$ for all $\bar{t}\in(\bT\cap\bar{f}\bK\bar{f}^{-1})^F$. Take
$t\in T^F$; as in the previous section, its eigenvalues
satisfy $\alpha_{(\xi,j,i+1)}=\alpha_{(\xi,j,i)}^q$. By \eqref{descripeqn},
$\bar{t}\in(\bT\cap\bar{f}\bK\bar{f}^{-1})^F$ if and only if, for some 
$\beta\in\Fq^\times$:
\begin{itemize}
\item for all $(\xi,j)\in\Lambda_{w_f}^1$, $\alpha_{(\xi,j,1)}^2=\beta$;
\item for all $(\xi,j)\in\Lambda_{w_f}^2$, 
$\alpha_{(\xi,j,1)}^{1+q^{\frac{1}{2}m_\xi(\nu_\xi)_j}}=\beta$;
\item for all $(\xi,j)\in\Lambda_{w_f}^3$,
$\alpha_{(\xi,j,1)}\alpha_{(\xi,\tilde{j},1)}^{q^{i(w_f,\xi,j)}}=\beta$.
\end{itemize}
If we let $\pm\sqrt{\beta}$ denote the two square roots of $\beta$ in 
$\Fqs^\times$, then
\beq
\begin{split} 
Z_{\bG}^\circ(\bar{t})\cap Z_{\bG}((\bT\cap\bar{f}\bK\bar{f}^{-1})^\circ)
&=(GL(\negthickspace\negthickspace\bigoplus_{\substack{(\xi,j,i)\\
(\xi,j)\in\Lambda_{w_f}^{1}\\\alpha_{(\xi,j,i)}=\sqrt{\beta}}} \negthickspace
L_{(\xi,j,i)})\times
GL(\negthickspace\negthickspace\bigoplus_{\substack{(\xi,j,i)\\
(\xi,j)\in\Lambda_{w_f}^{1}\\\alpha_{(\xi,j,i)}=-\sqrt{\beta}}} \negthickspace
L_{(\xi,j,i)})\\
&\qquad\times\prod_{\substack{(\xi,j,i)\\(\xi,j)\notin\Lambda_{w_f}^{1}}}
GL(L_{(\xi,j,i)}))/Z.
\end{split}
\eeq
If $\beta\in(\Fq^\times)^2$, then $(\pm\sqrt{\beta})^q=\pm\sqrt{\beta}$,
so the $\Fq$-rank of this group is the same as that of
$Z_{\bG}((\bT\cap\bar{f}\bK\bar{f}^{-1})^\circ)$. On the other hand, if 
$\beta\in\Fq^\times\setminus
(\Fq^\times)^2$, then $(\pm\sqrt{\beta})^q=\mp\sqrt{\beta}$, so $F$
interchanges the first two $GL$ factors, which means that the 
$\Fq$-rank of this group differs from that of
$Z_{\bG}((\bT\cap\bar{f}\bK\bar{f}^{-1})^\circ)$ by
$\frac{1}{2}\sum_{(\xi,j)\in\Lambda_{w_f}^{1}}m_\xi(\nu_\xi)_j$; moreover,
$m_\xi(\nu_\xi)_j$ must be even for all $(\xi,j)\in\Lambda_{w_f}^1$.
Thus
\beq  
\epsilon_{\bT,\bar{f}}(\bar{t})=
\left\{\begin{array}{cl}
1,&\text{ if $\beta\in(\Fq^\times)^2$,}\\
{\displaystyle (-1)^{\frac{1}{2}\sum_{(\xi,j)\in\Lambda_{w_f}^{1}}m_\xi(\nu_\xi)_j},}
&\text{ if $\beta\in\Fq^\times\setminus
(\Fq^\times)^2$,}
\end{array}\right.
\eeq
the second case being possible only when $w_f\in Y$.
Now consider $\lambda(\bar{t})=\lambda(t)$. If $\beta\in(\Fq^\times)^2$, then
\bes
\begin{split}
\lambda(t)
&=\prod_{(\xi,j)\in\Lambda_{w_f}^1}
\lng\alpha_{(\xi,j,1)}, N(\xi)\rng^{(\nu_\xi)_j}
\prod_{(\xi,j)\in\Lambda_{w_f}^2}
\lng\beta,N(\xi)\rng^{\frac{1}{2}(\nu_\xi)_j}
\prod_{(\xi,j)\in\widetilde{\Lambda_{w_f}^3}}
\lng\beta,N(\xi)\rng^{(\nu_\xi)_j}\\
&=\prod_{(\xi,j)\in\Lambda_{w_f}^1}
\lng\alpha_{(\xi,j,1)}\sqrt{\beta}^{-1}, N(\xi)\rng^{(\nu_\xi)_j},
\end{split}
\ees
where the first equality is by the reasoning before Proposition \ref{spbijprop},
and the second is by dividing by $\lng\sqrt{\beta},\Pi(\unnu)\rng=1$.
Since the elements $\alpha_{(\xi,j,1)}\sqrt{\beta}^{-1}$ for
$(\xi,j)\in\Lambda_{w_f}^1$ are all $\pm 1$,
and any choice of signs is possible, we see that $\lambda(\bar{t})=1$ for all such
$\bar{t}$ if and only if $w_f\in X$. So it only remains to consider the case when 
$w_f\in X\cap Y$
and $\bar{t}$ is such that $\beta\in\Fq^\times\setminus(\Fq^\times)^2$.
Since $\lng -1,\xi\rng_{m_\xi(\nu_\xi)_j}=1$ for all $(\xi,j)\in\Lambda_{w_f}^1$,
we have
\bes
\begin{split}
\lambda(t)&=\prod_{(\xi,j)\in\Lambda_{w_f}^1}
\lng\sqrt{\beta},\xi\rng_{m_\xi(\nu_\xi)_j}
\prod_{(\xi,j)\in\Lambda_{w_f}^2}
\lng\beta^{\frac{1}{2}(\nu_\xi)_j},\xi\rng_{m_\xi}
\prod_{(\xi,j)\in\widetilde{\Lambda_{w_f}^3}}
\lng\beta^{(\nu_\xi)_j},\xi\rng_{m_\xi}\\
&=\Phi(\unnu)\prod_{(\xi,j)\in\Lambda_{w_f}^2}
d_\xi^{\frac{1}{2}m_\xi(\nu_\xi)_j}
\prod_{(\xi,j)\in\widetilde{\Lambda_{w_f}^3}}
d_\xi^{m_\xi(\nu_\xi)_j},
\end{split}
\ees
where we have used the fact that $\lng\beta,\xi\rng_{m_\xi}=
d_\xi^{m_\xi}\lng\sqrt{\beta},\xi\rng_{2m_\xi}$. The result follows.
\epf

Combining \eqref{lusztigeqn}, Lemma \ref{partiallemma},
Proposition \ref{fibreprop}, and Proposition
\ref{mainprop}, we deduce that
\bes
\begin{split}
\langle B_\unnu, \Ind_{\bK^{F}}^{\bG^{F}}(1)\rangle
&=\frac{1}{2}(1+\epsilon\epsilon_{\unnu})\sum_{w\in W(T)_{\lambda,\ffi}^F}
\left\{\begin{array}{cl}
1&\text{ if $\Phi(w)=\Phi(\unnu)$}\\
0&\text{ if $\Phi(w)=-\Phi(\unnu)$}
\end{array}\right\}\\
&\quad+\sum_{w\in(X\cap Y)\setminus W(T)_{\lambda,\ffi}^F}
\negthickspace\negthickspace
2^{\ell_w^1-1}(-1)^{\ell_w^1}
\left\{\begin{array}{cl}
1&\text{ if $\Phi(w)=\Phi(\unnu)$}\\
0&\text{ if $\Phi(w)=-\Phi(\unnu)$}
\end{array}\right\}\\
&\quad+\sum_{w\in X\setminus Y} 2^{\ell_w^1-2}(-1)^{\ell_w^1}.
\end{split}
\ees
After some slight rearrangement, and noting that for
$w\in W(T)_{\lambda,\ffi}^F$ we have $\Phi(w)\Phi(\unnu)=\lng\beta,
\prod_{\xi\in\subsigL} N(\xi)^{\frac{1}{2}|\nu_\xi|}\rng$, this becomes
\beq \label{threetermeqn}
\begin{split}
\langle B_\unnu,& \Ind_{\bK^{F}}^{\bG^{F}}(1)\rangle
=\frac{1}{4}\sum_{w\in X}(-2)^{\ell_w^1}\\
&\quad+\left\{\begin{array}{cl}
\frac{1}{2}\epsilon\epsilon_{\unnu}\,|W(T)_{\lambda,\ffi}^F|,&\begin{array}{c}
\text{if all $|\nu_\xi|$ are even and}\\
\prod_{\xi\in\subsigL} N(\xi)^{\frac{1}{2}|\nu_\xi|}=1,\end{array}\\
0,&\text{ otherwise}
\end{array}\right.\\
&\quad+\frac{\Phi(\unnu)}{4}\sum_{w\in X\cap Y}\Phi(w)(-2)^{\ell_w^1}.
\end{split}
\eeq
The final step is to use the bijection
$W(T)_{\lambda,\inv}^F\isomto\prod_{\xi\in\subsigL} Z_{\inv}^{\nu_\xi}$
to identify the three terms of \eqref{threetermeqn} with the three terms
of Theorem \ref{pgothm}. 

The first has essentially already been done:
in \cite[\S4.1]{mysymm} it is observed that
\beq 
\sum_{w\in X}(-2)^{\ell_w^1}=\prod_{\substack{\xi\in\subsigL\\d_\xi=1}}
(\sum_{w_\xi\in Z_{\inv}^{\nu_\xi}}(-2)^{\ell_{w_\xi}^1(\nu_\xi)})
\prod_{\substack{\xi\in\subsigL\\d_\xi=-1}}
(\sum_{\substack{w_\xi\in Z_{\inv}^{\nu_\xi}\\
\ell_{w_\xi}^1(\nu_\xi)_1=0}}(-2)^{\ell_{w_\xi}^1(\nu_\xi)}),
\eeq
where $\ell_{w}^1(\nu)$ and $\ell_{w}^1(\nu)_1$ are defined as
in \cite[\S1.4]{mysymm}. Hence we need only invoke the combinatorial facts:
\begin{gather}
\label{combglnoneqn}
\sum_{w\in Z^{\nu}_\inv}
(-2)^{\ell_{w}^{1}(\nu)}
=(-1)^{|\nu|}\sum_{\rho\vdash|\nu|}
(\prod_i (m_i(\rho)+1))\chi_\nu^{\rho},\text{ and}
\\
\label{combglnsoneqn}
\sum_{\substack{w\in Z^{\nu}_\inv\\\ell_w^{1}(\nu)_1=0}}
\negthickspace\negthickspace\negthickspace
(-2)^{\ell_w^{1}(\nu)}
=\sum_{\substack{\rho\vdash|\nu|\\\rho' \text{ even}}}
\chi_\nu^{\rho}, 
\end{gather}
which are \cite[(4.1.2) and (4.1.3)]{mysymm}.

The second term of \eqref{threetermeqn} has also been done: in 
\cite[\S4.2]{mysymm} it is observed that
\beq 
\epsilon_{\unnu}\,|W(T)_{\lambda,\ffi}^F|=\prod_{\xi\in\subsigL}
\epsilon_{\nu_\xi}\,|Z_{\ffi}^{\nu_\xi}|,
\eeq
and using \eqref{macdonaldeqn} we obtain the second term of
Theorem \ref{pgothm}.

It only remains to analyse the third term of \eqref{threetermeqn}. It is clear
that under the bijection 
$W(T)_{\lambda,\inv}^F\isomto\prod_{\xi\in\subsigL} Z_{\inv}^{\nu_\xi}$, the
subset $X\cap Y$ of $W(T)_{\lambda,\inv}^F$ corresponds to the set
\[ \{(w_\xi)\in Z_{\inv}^{\nu_\xi}\,|\,\ell_{w_\xi}^1(\nu_\xi)_1=0\text{ whenever }
2\nmid m_\xi\text{ or }d_\xi=-1\}. \]
Therefore
\beq \label{lasteqn}
\begin{split}
\sum_{w\in X\cap Y}\Phi(w)(-2)^{\ell_w^1}&=
\prod_{\substack{\xi\in\langle\sigma\rangle\setminus L\\d_\xi=1\\2\nmid m_\xi}}
\sum_{\substack{w_\xi\in Z_{\inv}^{\nu_\xi}\\
\ell_{w_\xi}^1(\nu_\xi)_1=0}}
(-1)^{\ell_{w_\xi}^1(\nu_\xi)_{2\text{ mod }4}}
(-2)^{\ell_{w_\xi}^1(\nu_\xi)}\\
&\qquad\thickspace\thickspace\times
\prod_{\substack{\xi\in\langle\sigma\rangle\setminus L\\d_\xi=1\\2\mid m_\xi}}
(-1)^{\frac{1}{2}m_\xi|\nu_\xi|}
\sum_{w_\xi\in Z_{\inv}^{\nu_\xi}}
(-2)^{\ell_{w_\xi}^1(\nu_\xi)}\\
&\qquad\thickspace\thickspace\times
\prod_{\substack{\xi\in\langle\sigma\rangle\setminus L\\ d_\xi=-1}}
(-1)^{\frac{1}{2}m_\xi|\nu_\xi|}
\sum_{\substack{w_\xi\in Z_{\inv}^{\nu_\xi}\\
\ell_{w_\xi}^1(\nu_\xi)_1=0}}
(-2)^{\ell_{w_\xi}^1(\nu_\xi)}.
\end{split}
\eeq
In the first of the three groups of factors on the right-hand side, we have
used the fact that when $m_\xi$ is odd and $\ell_{w_\xi}^1(\nu_\xi)_1=0$,
\[ \frac{1}{2}\sum_{j\in\Lambda_{w_\xi}^1(\nu_\xi)}m_\xi(\nu_\xi)_j
\equiv 
\sum_{j\in\Lambda_{w_\xi}^1(\nu_\xi)}\frac{1}{2}(\nu_\xi)_j\equiv
\ell_{w_\xi}^1(\nu_\xi)_{2\text{ mod }4}\text{ (mod $2$),} \]
and in the second group of factors we have used the fact that when $m_\xi$ is even,
\[ \frac{1}{2}\sum_{j\in\Lambda_{w_\xi}^1(\nu_\xi)}m_\xi(\nu_\xi)_j
\equiv \frac{1}{2}m_\xi|\nu_\xi|\text{ (mod $2$).} \]
To finish the proof of Theorem \ref{pgothm},
we apply \eqref{combglnoneqn} to the second group of factors,
\eqref{combglnsoneqn} to the third group of factors, and the following identity
to the first group of factors:
\beq
\sum_{\substack{w\in Z^{\nu}_\inv\\\ell_w^{1}(\nu)_1=0}}
\negthickspace\negthickspace\negthickspace
(-1)^{\ell_{w}^1(\nu)_{2\text{ mod }4}}
(-2)^{\ell_w^{1}(\nu)}
=\sum_{\substack{\rho\vdash|\nu|\\2\mid m_{2i+1}(\rho), \forall i}}
(-1)^{\frac{1}{2}|\rho|+\ell(\rho)_{2\text{ mod }4}}
(\prod_i (m_{2i}(\rho)+1))\chi_{\nu}^{\rho}.
\eeq
This is merely \cite[(4.3.2)]{mysymm} with both sides multiplied by
$(-1)^{\frac{1}{2}|\nu|}$.


\begin{thebibliography}{10}

\bibitem{green}
{\sc J.~A. Green}, {\em The characters of the finite general linear groups},
  Trans. Amer. Math. Soc., 80 (1955), pp.~402--447.

\bibitem{mysymm}
{\sc A.~Henderson}, {\em Symmetric subgroup invariants in irreducible 
  representations of {$G^F$}, when {$G=GL_n$}}, J. Algebra, 261 (2003), 
  pp.~102--144. 

\bibitem{greenpolys}
{\sc G.~Lusztig}, {\em On the {Green} polynomials of classical groups},
  Proc. London Math. Soc., 33 (1976), pp.~443--475.

\bibitem{symmfinite}
\leavevmode\vrule height 2pt depth -1.6pt width 23pt, 
{\em Symmetric spaces
  over a finite field}, in The Grothendieck Festschrift, III, no.~88 in
  Progress in Mathematics, Birkhauser (Boston), 1990, pp.~57--81.

\bibitem{macdonald}
{\sc I.~G. Macdonald}, {\em Symmetric Functions and {Hall} Polynomials}, 
Oxford Univ. Press, second~ed., 1995.

\end{thebibliography}
\end{document}